\documentclass{ifacconf}
\usepackage{amsmath}  
\usepackage{amssymb}
\usepackage{upgreek}
\usepackage{dsfont}
\usepackage{algorithm}
\usepackage{algorithmic}
\usepackage{graphicx}      

\makeatletter
\let\old@ssect\@ssect 
\makeatother

\usepackage{natbib}
\usepackage{hyperref}

\makeatletter
\def\@ssect#1#2#3#4#5#6{%
  \NR@gettitle{#6}
  \old@ssect{#1}{#2}{#3}{#4}{#5}{#6}
}
\makeatother

\makeatletter
\def\endfigure{\end@float}
\def\endtable{\end@float}
\makeatother

\let\ifacconfcaptionwidth\captionwidth
\usepackage[caption=false]{subfig}
\let\captionwidth\ifacconfcaptionwidth

\DeclareMathOperator{\sgn}{sgn}
\DeclareMathOperator{\Tr}{Tr}
\DeclareMathOperator{\prox}{prox}

\begin{document}
\begin{frontmatter}

\title{Sparse optimal control of networks with multiplicative noise via policy gradient\thanksref{footnoteinfo}} 

\thanks[footnoteinfo]{This material is based on work supported by the Army Research Office under grant W911NF-17-1-0058 and the National Science Foundation under grant CMMI-1728605.}

\author{Benjamin Gravell} 
\author{\quad\quad Yi Guo} 
\author{\quad\quad Tyler Summers}

\address{The University of Texas at Dallas, 
   Richardson, TX 75080 USA 
   }

\begin{abstract}                
We give algorithms for designing near-optimal sparse controllers using policy gradient with applications to control of systems corrupted by multiplicative noise, which is increasingly important in emerging complex dynamical networks. Various regularization schemes are examined and incorporated into the optimization by the use of gradient, subgradient, and proximal gradient methods. Numerical experiments on a large networked system show that the algorithms converge to performant sparse mean-square stabilizing controllers.
\end{abstract}


\end{frontmatter}

\section{Introduction}
Emerging highly distributed networked dynamical systems, such as critical infrastructure for power, water, and transportation, are high-dimensional and increasingly instrumented with new sensing, actuation, and communication technologies. A key problem is to design high performance control architectures that limit the number of actuators, sensors, and actuator-sensor communication links to reduce complexity and cost. Sparse control architectures may be crucial for managing complexity in emerging complex networks, but require solution of extremely difficult mixed combinatorial-continuous optimization problems.

There is a variety of performance metrics and optimization methodology for sparse control architecture design in the recent literature. Examples include structural rank conditions from \cite{liu2011controllability,ruths2014control,olshevsky2014minimal}, controllability and observability Gramians from \cite{Pasqualetti2014c,summers2014submodularity,tzoumas2016,jadbabaie2018deterministic}, and optimal and robust control metrics from \cite{Hassibi1998,Polyak-LMI_sparse_fb,jovanovic2016controller,summers2016actuator,Taha2017d,Zare2018CDC}, which are optimized via greedy algorithms, convex and mixed-integer optimization, and randomization.

Here we develop methods for sparse optimal control design in dynamical networks with multiplicative noise via policy gradient algorithms with sparsity-inducing regularization.  Multiplicative noise arises in many networked systems when the weights of edges connecting nodes are stochastic in time. 
The noise is thus on the system parameters themselves and has a fundamentally different effect on the state evolution than additive noise, and indeed can lead to dramatic robustness issues. Specifically, 
a noise-ignorant classical optimal linear-quadratic (LQ) controller may actually \textit{destabilize} a multiplicative noise system in the mean-square sense, even if the system was open-loop mean-square stable. Therefore noise-aware control is imperative to network performance and robustness. Moreover, the policy gradient methods we propose here, which operate directly on policy parameters, facilitate data-driven sparse control design when the model is unknown, a topic we are exploring in ongoing work.




In Section 2 we formulate the problem and discusses a policy gradient approach to optimal control design for linear-quadratic systems with multiplicative noise. In Section 3 we propose several sparse control design methods for sensor and actuator selection and communication network design using gradient, subgradient, and proximal algorithms. In Section 4 we present numerical experiments to illustrate the results. Section 5 concludes.


\section{Problem formulation}
Consider the discrete-time linear quadratic regulator with multiplicative noise (LQRm) optimal control problem
\begin{equation} \label{eq:LQRm}
\begin{aligned}
&\underset{{\pi \in \Pi}}{\text{min}} && J(\pi) = \mathds{E}\sum_{t=0}^\infty (x_t^T Q x_t + u_t^T R u_t), \\
&\text{s.t.} && x_{t+1} = (A  + \sum_{i=1}^p \delta_{it} A_i) x_t +  (B + \sum_{j=1}^q \gamma_{jt} B_j) u_t
\end{aligned}
\end{equation}
where $x_t \in \mathds{R}^n$ is the system state, $u_t \in \mathds{R}^m$ is the control input, $x_0$ is randomly distributed according to $\mathcal{P}$, expectation is with respect to $x_0,\delta_{it}, \gamma_{jt}$, and $Q\succeq 0$ and $R\succ 0$. The dynamics incorporate multiplicative noise terms modeled by the mutually independent and i.i.d. (over time) zero-mean random variables $\delta_{it}$ and $\gamma_{jt}$, which have variance $\alpha_i$ and $\beta_j$, respectively. The matrices $A_i \in \mathds{R}^{n \times n}$ and $B_i \in \mathds{R}^{n \times m}$ specify how each noise term affects the system dynamics and input matrices. The goal is to determine an optimal closed-loop feedback policy $\pi$  with $u_t = \pi(x_t)$. We assume that the problem data $A$, $B$, $\alpha_i$, $A_i$, $\beta_j$, and $B_j$ are such that the optimal value of the problem exists and is finite. Feasibility of this problem is ensured if the system is mean-square stabilizable.
\begin{defn}[Mean-square stability]
    The system in \eqref{eq:LQRm} is stable in the mean-square sense if  
        $\lim_{t \to \infty}\mathds{E}[x_t x_t^T] = 0$
    for any given initial covariance $\mathds{E}x_0 x_0^T$.
\end{defn}

We are ultimately interested in the problem
\begin{equation} \label{eq:LQRm_sparse}
\begin{aligned}
&\underset{{\pi \in \Pi}}{\text{min}} && J(\pi) + \gamma J_{\text{reg}}(\pi)
\end{aligned}
\end{equation}
where $J_{\text{reg}}(\pi)$ is a sparsity-promoting regularizer of the policy $\pi$ and $\gamma$ specifies the importance of sparsity. The regularizer ideally would measure the number of actuators, sensors, or actuator-sensor links, but for computational tractability will be replaced by other functions defined later. We begin by discussing the solution for $\gamma=0$.


\subsection{Optimal control via value iteration}
Dynamic programming can be used to show that the optimal policy is linear state feedback with $u_t = K^* x_t$
where $K^* \in \mathds{R}^{m \times n}$ and the resulting optimal cost for a fixed initial state is quadratic with $V_{K^*}(x_0) = x_0^T P x_0$ where $P \in \mathds{R}^{n \times n}$ is a symmetric positive definite matrix. When the model parameters are known, there are several known ways to compute the optimal feedback gains and corresponding optimal cost. The optimal cost is given by the solution of the generalized algebraic Riccati equation (ARE) (see, e.g., \cite{Damm2004}).  
\begin{align} \label{genriccati}
P &= Q + A^T P A + \sum_{i=1}^p \alpha_i A_i^T P A_i \\
&\quad -  A^T P B (R + B^T P B + \sum_{j=1}^q \beta_j B_j^T P B_j)^{-1} B^T P A. \nonumber 
\end{align}
This can be solved via value iteration, 
and the optimal gain matrix is
\begin{equation}
    K^* = -  \Big(R + B^T P B + \sum_{j=1}^q \beta_j B_j^T P B_j \Big)^{-1} B^T P A.
\end{equation}

\subsection{Optimal control via policy gradient}

For a fixed mean-square stabilizing linear state feedback policy $u_t = K x_t$, there exists a positive semidefinite cost matrix $P_K$ which characterizes the cost by
\begin{equation}
    J(K) = \underset{x_0}{\mathds{E}} x_0^T P_{K} x_0
\end{equation}
and is the solution to the generalized Lyapunov equation
\begin{align}
    &P_{K} =  Q + K^TRK  + (A+BK)^T P_{K} (A+BK) \nonumber \\
    &\quad + \sum_{i=1}^p \alpha_i A_i^T P_K A_i + \sum_{j=1}^q \beta_j K^T B_j^T P_K B_j K. \label{eq:gen_lyap_P}
\end{align}
Furthermore, there exists a positive semidefinite infinite-horizon aggregate state covariance matrix
\begin{align*}
\Sigma_{K} = \underset{x_0,\delta_{it},\gamma_i}{\mathds{E}} \sum_{t=0}^\infty x_t x_t^T
\end{align*}
which is the solution to the generalized Lyapunov equation
\begin{align}
    \Sigma_{K} &= \Sigma_0 + (A+BK) \Sigma_{K} (A+BK)^T + \sum_{i=1}^p \alpha_i A_i \Sigma_K A_i^T \nonumber \\
    &\quad + \sum_{j=1}^q \beta_j B_j K \Sigma_K K^T B_j^T, \label{eq:gen_lyap_S}
\end{align}
where $\Sigma_{0} = \underset{x_0}{\mathds{E}}\left[x_0x_0^T\right]$. 
Thus, we have
\begin{align}
    J(K) = \Tr((Q+K^T R K) \Sigma_K) = \Tr(P_K \Sigma_0). \label{eq:cost}
\end{align}


This leads to the idea of performing gradient descent on $J$ (i.e., policy gradient) to find the optimal gain matrix:
\begin{equation}
    K \leftarrow K - \eta \nabla_K J(K),
\end{equation}
for a fixed step size $\eta$.
In this work we consider only the case where the model parameters are known, but the methods presented are immediately usable in the model-unknown case by estimating the gradient from trajectory data. The policy gradient for linear state feedback policies applied to the LQRm problem has the following form:
\begin{lem} The LQRm policy gradient is given by
\begin{align}
&\nabla_K J(K) \\
&= 2\bigg[ \Big(R + B^T P_{K} B + \sum_{j=1}^q \beta_j B_j^T P_K B_j \Big) K + B^T P_{K} A\bigg] \Sigma_{K} . \nonumber
\end{align}
\end{lem}
The proof is omitted due to space constraints and can be found in our technical report (see \cite{Gravell2019unpublished}).

\subsection{Gradient domination}
It was shown recently by \cite{Fazel2018} that although the deterministic LQR cost is nonconvex, it is \emph{gradient dominated}, also known as the Polyak-{\L}ojasiewicz inequality originally due to \cite{Polyak1963}.
It is simple to show that if a function has a Lipschitz continuous gradient and satisfies this condition then performing gradient descent with a sufficiently small constant step size will result in asymptotic convergence to the optimal function value at a linear rate (see \cite{Karimi2016}). For the LQRm problem, so long as the initial controller is stabilizing the LQRm cost is continuously differentiable over the sublevel set associated with the initial controller and thus the gradient possesses a local Lipschitz constant $L$ on this set. Identifying $L$ and the gradient domination constant is necessary for selection of a step size which is guaranteed to give convergence using gradient descent. Quantifying these constants is difficult but possible via lengthy chains of matrix inequalities as demonstrated by \cite{Fazel2018}. 

These results extend readily to the LQRm problem with relevant quantities pertaining to Lipschitz continuity of the gradient and the gradient domination conditions modified suitably to accommodate the multiplicative noise. In particular, the effect of the noise is to decrease the maximum step size that can be taken using gradient descent. We now state the relevant lemmas; the proofs are lengthy and can be found in our technical report (see \cite{Gravell2019unpublished}).

\begin{lem}[Gradient domination] \label{lemma:gradient_dominated} 
    The LQRm cost $J(K)$ satisfies the gradient domination condition
    \begin{align}
        J(K) - J(K^*) &\leq \frac{\|\Sigma_{K^*}\|}{4\sigma_{\text{min}}(R) \sigma_{\text{min}}(\Sigma_{0})^2 } \|\nabla_K J(K)\|_F^2
    \end{align}     
\end{lem}

\begin{lem}[Gradient descent, convergence rate] \label{lemma:grad_exact_convergence}
    \ \\ Using the policy gradient step update
    \begin{equation}
        K^{(k+1)} = K^{(k)} - \eta \nabla_K J(K^{(k)})
    \end{equation}
    with step size $0 < \eta \leq c_{pg}$ gives global convergence to the optimal $K^*$ at a linear rate described by
    \begin{equation}
      \frac{ J(K^{(k+1)}) - J(K^*)}{ J(K^{(k)}) - J(K^*)} \leq 1 - \eta \frac{ \sigma_{\text{min}}(R)\sigma_{\text{min}}(\Sigma_{0})^2}{\|\Sigma_{K^*}\|} 
    \end{equation}
    where $c_{\text{pg}}$ is a constant which is polynomial in the parameters $A$, $B$, $B_j$, $Q$, $R$, $J(K^{(0)})$.
    
\end{lem}




\section{Sparse control design}
Entrywise, row, and column sparsity in $K$ correspond to actuator-sensor communication, actuator, and sensor sparsity respectively. With this in mind, we seek to solve the optimization problem of finding the sparsest set of entries, rows and/or columns of $K$ that achieve some prescribed level of performance in terms of the LQRm cost. However this problem is a nonconvex combinatorial problem which is NP-hard; the number of independent problem instances which must be solved scales factorially with $n$ and/or $m$. We instead turn to regularization as a heuristic to identifying good sparsity patterns.

\subsection{Insufficiency of naive hard thresholding}
The most na\"ive method of inducing sparsity is hard thresholding of the ARE solution as $K_{ij} = 0 \text{ if } |K_{ij}| < r$. However, in general this is not useful since the resulting gains may not be stabilizing. Consider the following example system without multiplicative noise:
\begin{align*}
    A = 
    \begin{bmatrix}
    0.4 &  0.9 & -0.3 \\
    0.7 & -0.3 & -0.4 \\
    -0.2 & 0.1 & -0.8
    \end{bmatrix}, \quad
    B = 
    \begin{bmatrix}
    0.2 & -0.6 \\
    -1.3 & -1.6 \\ 
    -0.3 & -1.5
    \end{bmatrix}, \\
    \mathcal{P} = \mathcal{N}(0,I_3),\quad Q=I_3,\quad R=I_2
\end{align*}
where $I_n$ is an $n \times n$ identity matrix. Imposing a hard threshold of $0.4$ on the ARE solution results in
\begin{align*}
    K = 
    \begin{bmatrix}
    0.503931 & -0.880322  & 0          \\
    0         &  0.614382 & -0.677758
    \end{bmatrix}
\end{align*}
which gives a closed-loop state transition matrix $A+BK$ with an eigenvalue of $1.048223$ outside the unit circle.
By contrast, by working with the regularized LQRm cost the optimal gains are always guaranteed to be stabilizing; even in the limit as the regularization weight $\rightarrow \infty$ the sparsity increases until the sparsest stabilizing solution is obtained. In practice, using a small step size helps ensure that each iterate remains inside the domain of $J(K)$.

\subsection{Regularization}
Certain types of regularization are well-known to be capable of inducing sparsity in the solutions to optimization problems. Perhaps the most basic and well-known is $l_1$-norm regularization which operates on a vector of decision variables; see \cite{Tibshirani1996} for the seminal LASSO problem for sparse least-squares model selection and \cite{Hassibi1998} for sparse control design. In the case of a convex objective, increasing the regularization weight tends to increase sparsity by moving the global minimum onto the coordinate axes. Once the regularized problem has been solved, a sparsity pattern can easily be identified from the (near-)zero entries.
In the current work we consider only the problem of identifying sparsity patterns, however an additional ``polishing'' step which involves re-solving the LQRm problem under the sparsity pattern can be performed to further improve the LQRm cost, as in \cite{Lin2013}. 

Entrywise sparsity is induced by the vector $l_1$-norm
\begin{align}
    \|K\|_{1} = \sum_{i=1}^m \sum_{j=1}^n |K_{ij}|.
\end{align}
Row and column sparsity are induced by using matrix row and column norms respectively defined as
\begin{align}
    \|K\|_{r} = \sum_{i=1}^m \|K^{r,i}\|_{\infty}, \ \text{ and } \  \|K\|_{c} = \sum_{i=1}^n \|K^{c,i}\|_{\infty},
\end{align}
where $\|K^{r,i}\|_\infty$ and $\|K^{c,i}\|_\infty$ are the maximum absolute values of the $i^{th}$ row and column respectively of $K$.
Row and column sparsity are also induced by the row and column group LASSO
\begin{align}
    \|K\|_{glr} = \sum_{i=1}^m \|K^{r,i}\|_{2}, \ \text{ and } \  \|K\|_{glc} = \sum_{i=1}^n \|K^{c,i}\|_{2},
\end{align}
where $\|K^{r,i}\|_\infty$ and $\|K^{c,i}\|_\infty$ are the vector $l_2$-norms of the $i^{th}$ row and column respectively of $K$.
Combined row and column sparsity can be induced by the row and column sparse group LASSO 
\begin{align}
    \|K\|_{sglr} = (1-\mu) \|K\|_{1} + \mu \|K\|_{glr} \\
    \|K\|_{sglc} = (1-\mu) \|K\|_{1} + \mu \|K\|_{glc}
\end{align}
or by various other weighted combinations of entrywise, row, and column norms.
We refer to $\|K\|_M$ as a generic nondifferentiable sparsity-inducing regularizer.

\subsection{Stationary point characterization}
Before proceeding, we must point out an important consequence of regularizing the LQRm cost. The sum of a convex function and a gradient dominated function is not gradient dominated in general, and in fact can have multiple local minima. For example, consider the scalar function
\begin{align}
    f(x) = x^2 + 4 ((x-8)^2 + 3 \sin^2(x-8))
\end{align}
where $x^2$ is strongly convex and $4((x-8)^2 + 3 \sin^2(x-8))$ is gradient dominated. But $f(x)$ has two local minima at $x=5.372$ and $x=7.459$ and therefore is not gradient dominated.

As a result any local first-order search procedure, such as those used by our algorithms, will not be guaranteed to find the global minimum. We conjecture that for the regularized LQRm problem there are at most two local minima, one associated with the LQRm cost and one associated with the regularization which tends to be more sparse. If this is so then choosing the initial point carefully may help the local search find the desired (sparser) local minimum. For open-loop mean-square systems, this motivates using zero gains as the initial condition. Likewise, in both the open-loop mean-square stable and unstable cases, an effective heuristic is to use the solution to a highly regularized problem instance to ``warm start'' another nearby problem instance with reduced regularization weight.

\subsection{Other step directions}
Promising choices of step directions other than the gradient $\nabla_K J(K)$ are the natural gradient $\nabla_K J(K)\Sigma_{K}^{-1}$ and the Gauss-Newton step $R_{K}^{-1} \nabla_K J(K)\Sigma_{K}^{-1}$ as given by \cite{Fazel2018}. When $\gamma=0$, these step directions give faster convergence than the gradient step and in fact the convergence proofs are much simpler than that for the gradient step. Unfortunately, adding a regularizer makes these steps more difficult to calculate; it is not simply the sum of the gradient of the regularizer and the unregularized natural gradient or Gauss-Newton step of the LQRm cost. For this reason we restrict our attention to the standard (sub)gradient directions.

\subsection{Regularized policy subgradient descent}
In order to use nondifferentiable regularizers we use subgradient methods which take steps in the direction of subgradients. It is known that using a constant step size gives convergence to a bounded neighborhood of the optimum and that a diminishing step size gives asymptotic, albeit slow, convergence (see \cite{Nesterov2013}).
One immediate issue is that subgradients are defined only for convex functions; since the LQRm cost is nonconvex, subgradients do not exist for the regularized LQRm cost. However we simply use the gradient of the LQRm cost plus the subgradient of the regularizer as the step direction. Thus our subgradient descent update is 
\begin{alg}[Policy subgradient update]
    \begin{align*}
        K^{(k+1)} &= K^{(k)}-\eta (\nabla_K J(K^{(k)}) + \gamma \blacktriangledown_K  \|K^{(k)}\|_M) \\
        K_{\text{best}} &= \text{argmax } \{C(K^{(k+1)}),C(K_{\text{best}})\}
    \end{align*}
    where $C(K)=J(K) + \gamma \|K\|_M$ and $\blacktriangledown_K$ is a subgradient.
\end{alg}
Another issue is that there is no guarantee of feasibility of each next step; it is possible to take a step so large that the next point is a mean-square unstable controller giving infinite objective cost. It is not straightforward to obtain restrictions on the step size to guarantee this feasibility. Gradient descent does not suffer from this problem since the gradient is guaranteed to be a true descent direction so there is always a sufficiently small step size to give a feasible next step. Nevertheless, in practice it is rare for a sufficiently small subgradient step to be infeasible. 

\subsection{Proximal policy gradient}
Proximal gradient methods have become a preferred way to solve optimization problems of the form
\begin{align*}
    \underset{x}{\text{min}} \ f(x) + g(x)
\end{align*}
where $f(x)$ has a Lipschitz continuous gradient and $g(x)$ is convex and nondifferentiable, as is the case when $g(x)$ is a sparsity-inducing regularizer. The proximal gradient method update is
\begin{align*}
    x^{(k+1)} =\prox_{\eta g}\Big(x^{(k)}-\eta \nabla f\big(x^{(k)}\big)\Big)
\end{align*}
where the proximity operator is defined as
\begin{align*}
    \prox_{\eta g}(v)=\underset{x}{\operatorname{argmin}}\left\{g(x)+(1 /( 2 \eta))\|x-v\|_{2}^{2}\right\} .
\end{align*}
Much of the existing literature examines the case where $f(x)$ is convex, in which case gradient descent is guaranteed to converge. The proximal operator has closed-form expressions for $\|K\|_{1}$ and $\|K\|_{glr}$ called \textit{soft thresholding} and \textit{block soft thresholding} (see \cite{Parikh2014}).
Thus to solve \eqref{eq:LQRm_sparse} we also use a proximal policy gradient algorithm:
\begin{alg}[Proximal policy gradient update]
    \begin{align*}
        K_s^{(k)} &= K^{(k)}-\eta \Delta K^{(k)} \\
        K^{(k+1)} &=\prox_{\eta \|K\|_M}\big(K_s^{(k)}\big)
    \end{align*}
    where $\Delta K^{(k)}$ is a generic step direction.
\end{alg}
A result from \cite{Hassan2018} guarantees convergence at a linear rate to the optimal function value using the proximal gradient method on a function satisfying a proximal gradient domination condition. This condition was shown to be equivalent by one given by \cite{Karimi2016} and an inequality from \cite{Kurdyka1998}. However, this condition is not guaranteed to hold when $f(x)$ is gradient dominated and $g(x)$ is convex; the full condition must be checked, which involves interaction between $f(x)$ and $g(x)$. It is nontrivial to verify that the condition is satisfied for the regularized LQRm cost. Empirically it appears that the inequality may be satisfied since the proximal gradient method converged to solutions similar to those from our other two methods.

\subsection{Regularized policy gradient descent}
Another algorithm for solving \eqref{eq:LQRm_sparse} is gradient descent:
\begin{alg}[Policy gradient update]
    \begin{align*}
        K^{(k+1)} &= K^{(k)}-\eta \nabla_K (J(K^{(k)}) + \gamma \|K\|_M)
    \end{align*}
\end{alg}
Here we use differentiable Huber-type losses $\|K\|_{M,h,\phi}$ in place of nondifferentiable regularizers, which replace linear corners with quadratic tips for decision variable values smaller than a specified threshold. Although the solutions produced are not exactly sparse, in practice entries are sufficiently close to zero to identify the sparsity pattern. Furthermore, by iteratively decreasing the threshold the solutions can be made arbitrarily close to truly sparse.

We define the Huber function of a scalar $a$ as
\begin{align*}
    h_{\phi}(a) = \left\{\begin{array}{ll} {|a|-\frac{1}{2} \phi} & {\text { if } |a| > \phi} \\ {\frac{1}{2 \phi} a^{2}} & {\text { if } |a| \leq \phi}\end{array}\right.
\end{align*}
and the $p$-Huber function (like a $p$-norm) of a vector $b$ as
\begin{align*}
    h_{p,\phi}(b) = \left\{\begin{array}{ll} {\|b\|_p-\frac{1}{2} \phi} & {\text { if } \|b\|_p > \phi} \\ {\frac{1}{2 \phi} \|b\|_p^{2}} & {\text { if } \|b\|_p \leq \phi}\end{array}\right.
\end{align*}
We define the vector Huber loss as 
\begin{align*}
    \|K\|_{1,h,\phi} = \sum_{i=1}^m \sum_{j=1}^n h_{\phi}(K_{ij})
\end{align*}
the Huber row and column norms as 
\begin{align*}
    \|K\|_{r,h,\phi} = \sum_{i=1}^m h_{\infty,\phi}(K^{r,i}) , \quad \|K\|_{c,h,\phi} = \sum_{i=1}^n h_{\infty,\phi}(K^{c,i})
\end{align*}
and the Huber row and column group LASSO as 
\begin{align*}
    \|K\|_{glr,h,\phi} = \sum_{i=1}^m h_{2,\phi}(K^{r,i}) , \quad \|K\|_{glc,h,\phi} = \sum_{i=1}^n h_{2,\phi}(K^{c,i})
\end{align*}
Subgradients of two regularizers and the gradients of their differentiable counterparts are given in Table~\ref{table:1}.
\begin{table}[ht!]
\begin{center}
\caption{Regularizer (sub)gradients}\label{table:1}
 \begin{tabular}{|c | c | c|} 
 \hline 
 \rule{0pt}{2ex} $\|K\|_M$ & $\blacktriangledown_K \|K\|_M$ & $\nabla_K \|K\|_{M,h,\phi}$\\ [0.5ex] 
 \hline\hline
 \rule{0pt}{2.5ex} $\|K\|_{1,ij}$ & $\sgn(K_{ij})$ & $K_{ij}/\max(|K_{ij}|,\phi)$ \\ [0.5ex] 
 \hline
 \rule{0pt}{2.5ex} $\|K\|_{glr}^{r,i}$ & $K^{r,i}/\|K^{r,i}\|_2$ &  $K^{r,i}/\max(\|K^{r,i}\|_2, \phi)$ \\ [0.5ex] 
 \hline
\end{tabular}
\end{center}
\end{table}

\section{Simulation results}
We considered an example system which represents diffusion dynamics on a particular undirected Erd\H os-R\'enyi random graph. It is well known that if $p_{ER}=(\log n+c) / n$ for constant $c \in \mathbb{R}$, then $ \lim_{n \rightarrow \infty} P(G(n, p) \text{ connected}) = e^{-e^{-c}}$ so we chose $n=51$, $c=7$ and $p_{ER}=0.2144$ and with probability $P=0.999$ obtained a connected graph (see \cite{Bollobas2001}). The graph was selected so that it was connected, which ensured controllability. The first row and and column of the graph Laplacian were removed in order to fix the system's state reference to the first node which removed the zero eigenvalue otherwise present. The continuous time system was discretized using a standard bilinear transform (Tustin's approximation) which preserves the open-loop mean stability of this system. Two multiplicative noises act each on $A$ and $B$ whose entries were drawn from a Gaussian distribution. The multiplicative noise variances were set at two levels, low and high, so that the system was open-loop mean-square stable and unstable, respectively.

\begin{figure}[ht!]
  \centering
  \subfloat[$75.5\%$ sparsity, $\gamma {=} 10$]{\includegraphics[width=4.3cm]{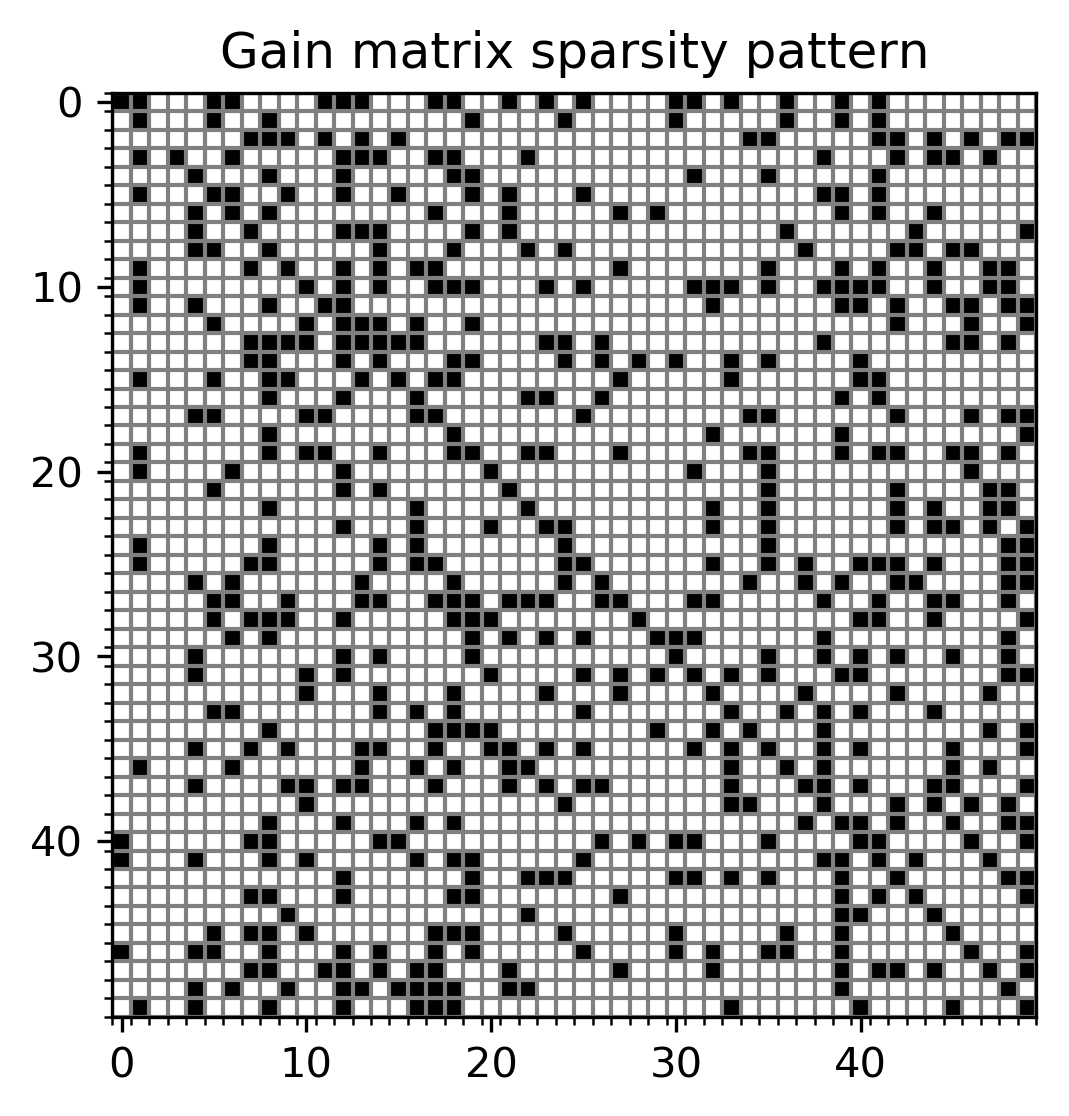}}
  \subfloat[$94.3\%$ sparsity, $\gamma {=} 320$]{\includegraphics[width=4.3cm]{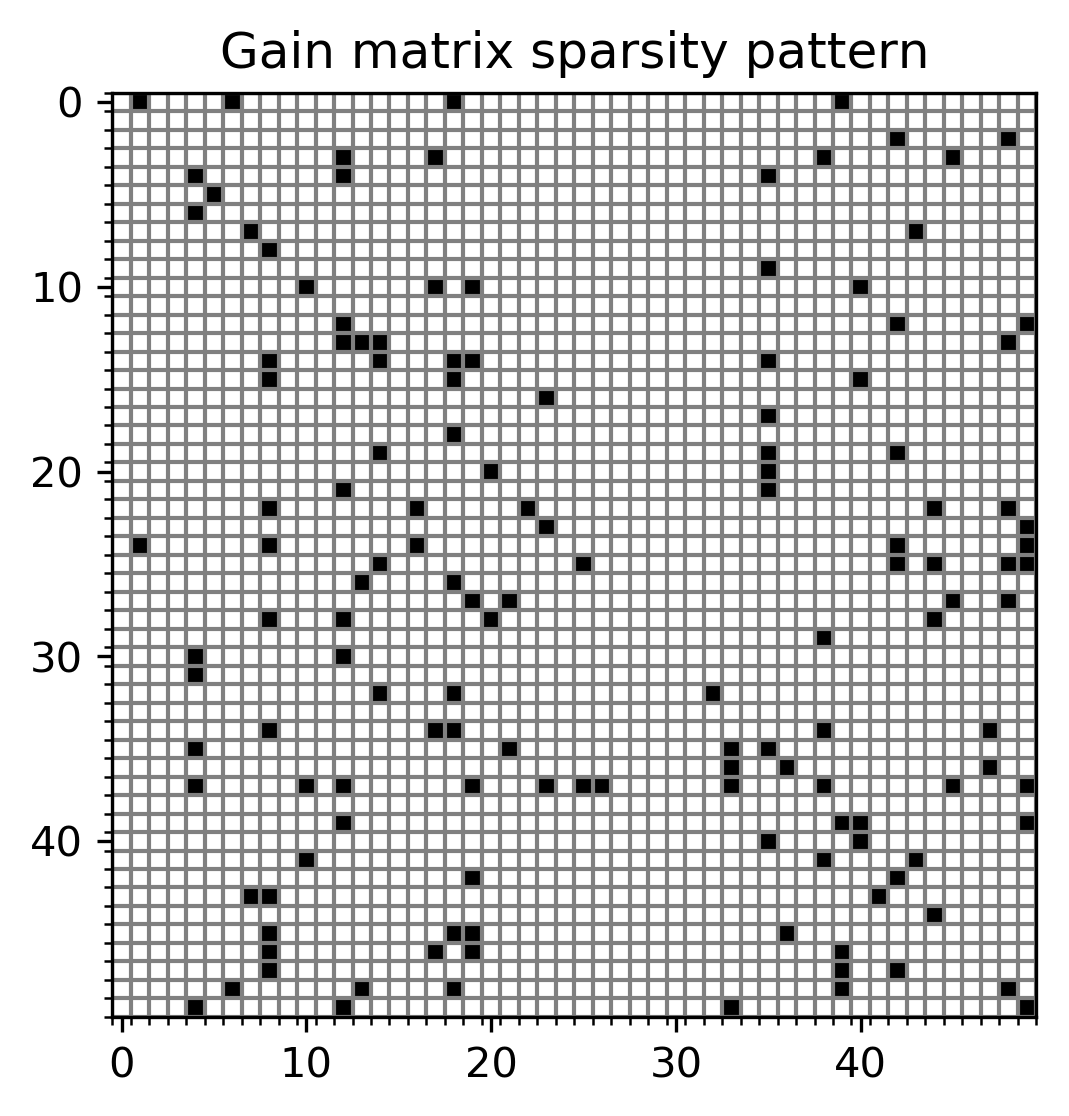}}
  \caption{Sparsity patterns for low noise, subgradient descent on the $l_1$-norm regularized LQRm cost.}
  \label{fig:sparsity_vec1}
\end{figure}
\begin{figure}[ht!]
    \centering
    \subfloat[$18\%$ sparsity, $\gamma {=} 400$]{\includegraphics[width=4.3cm]{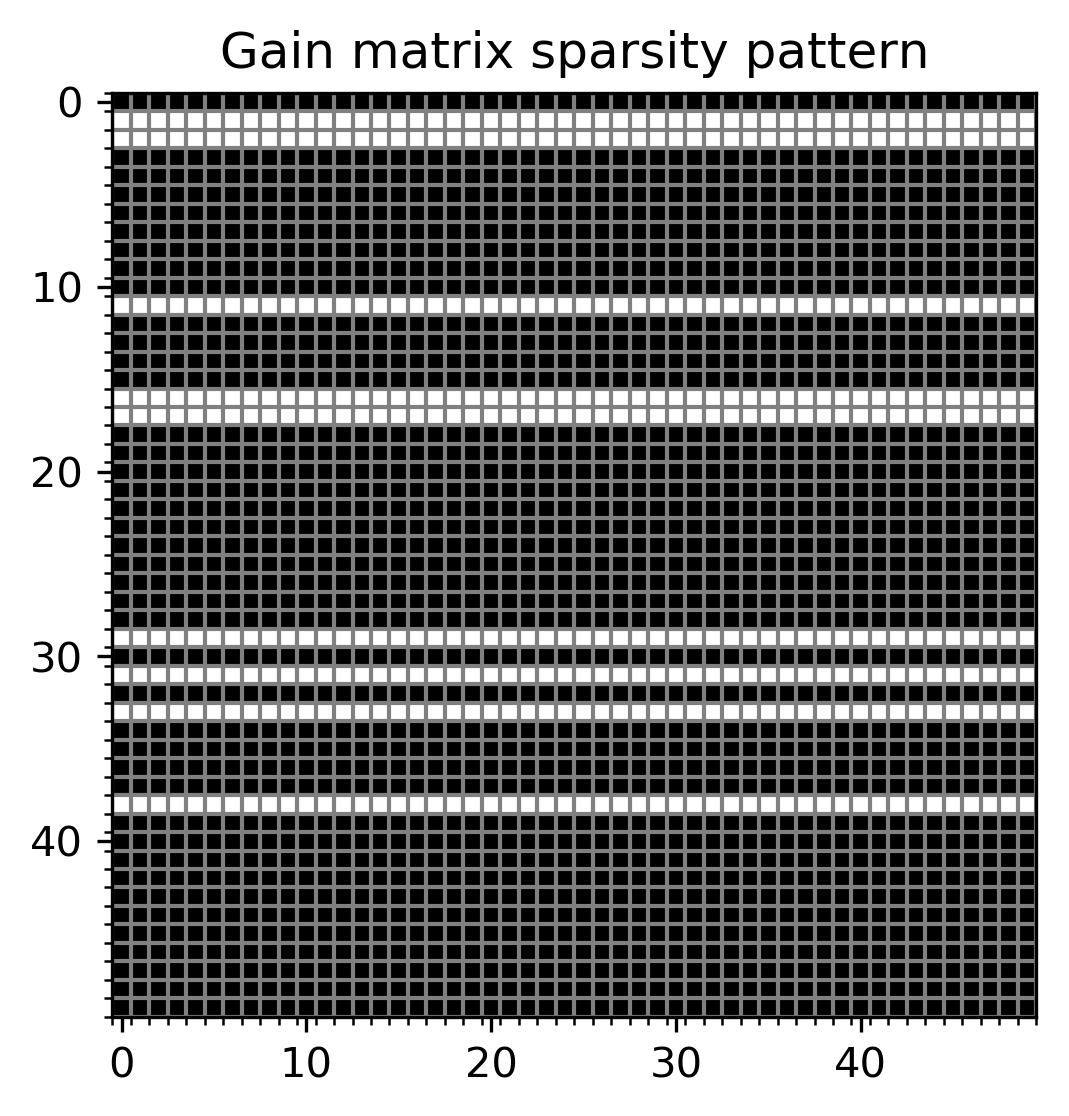}}
    \subfloat[$94\%$ sparsity, $\gamma {=} 18102$]{\includegraphics[width=4.3cm]{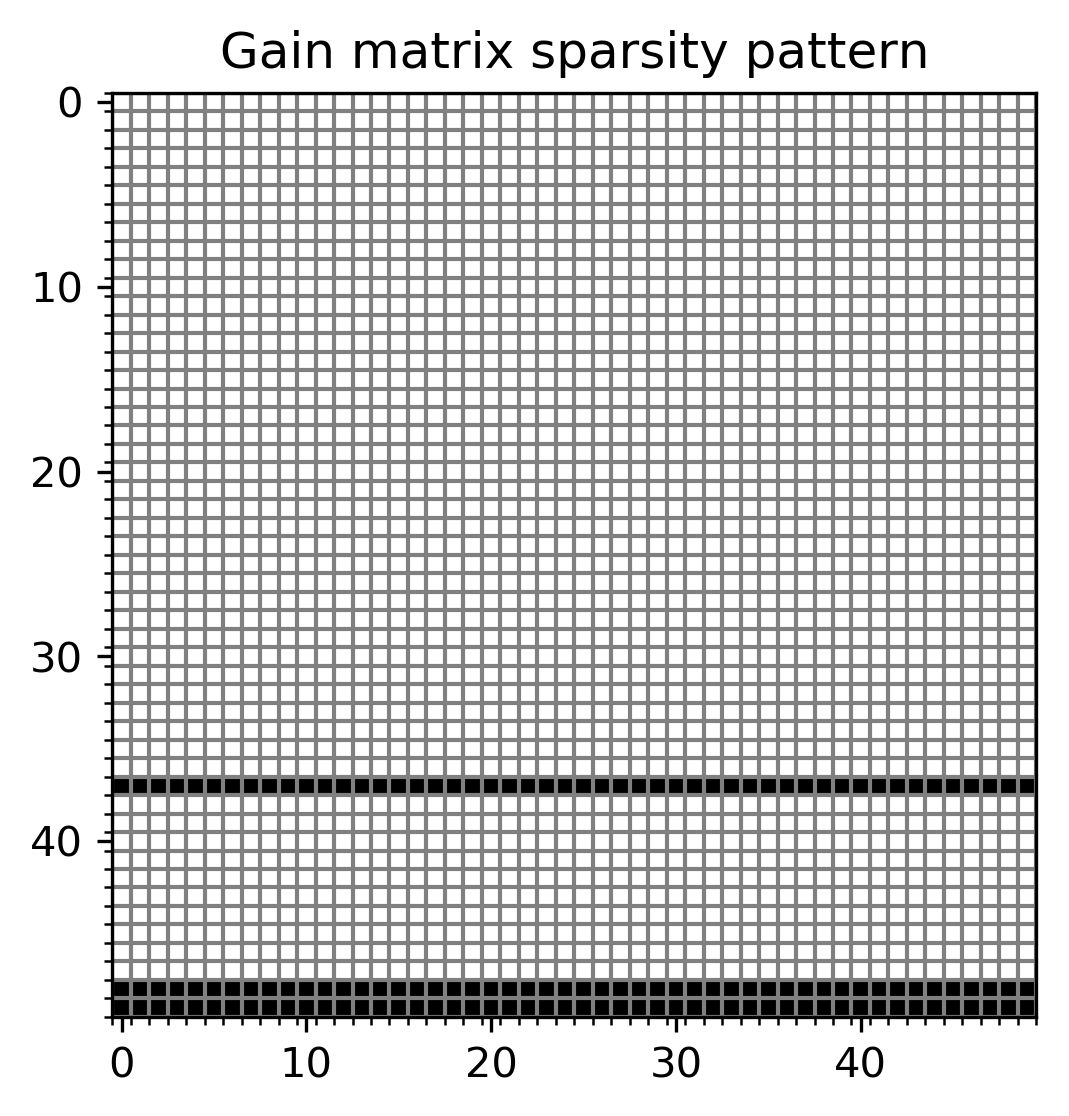}}
    \caption{Sparsity patterns for low noise, subgradient descent on the row group LASSO regularized LQRm cost.}
    \label{fig:sparsity_glr}
\end{figure}
\begin{figure}[ht!]
    \centering
    \subfloat[Low noise setting]{\includegraphics[width=4.3cm]{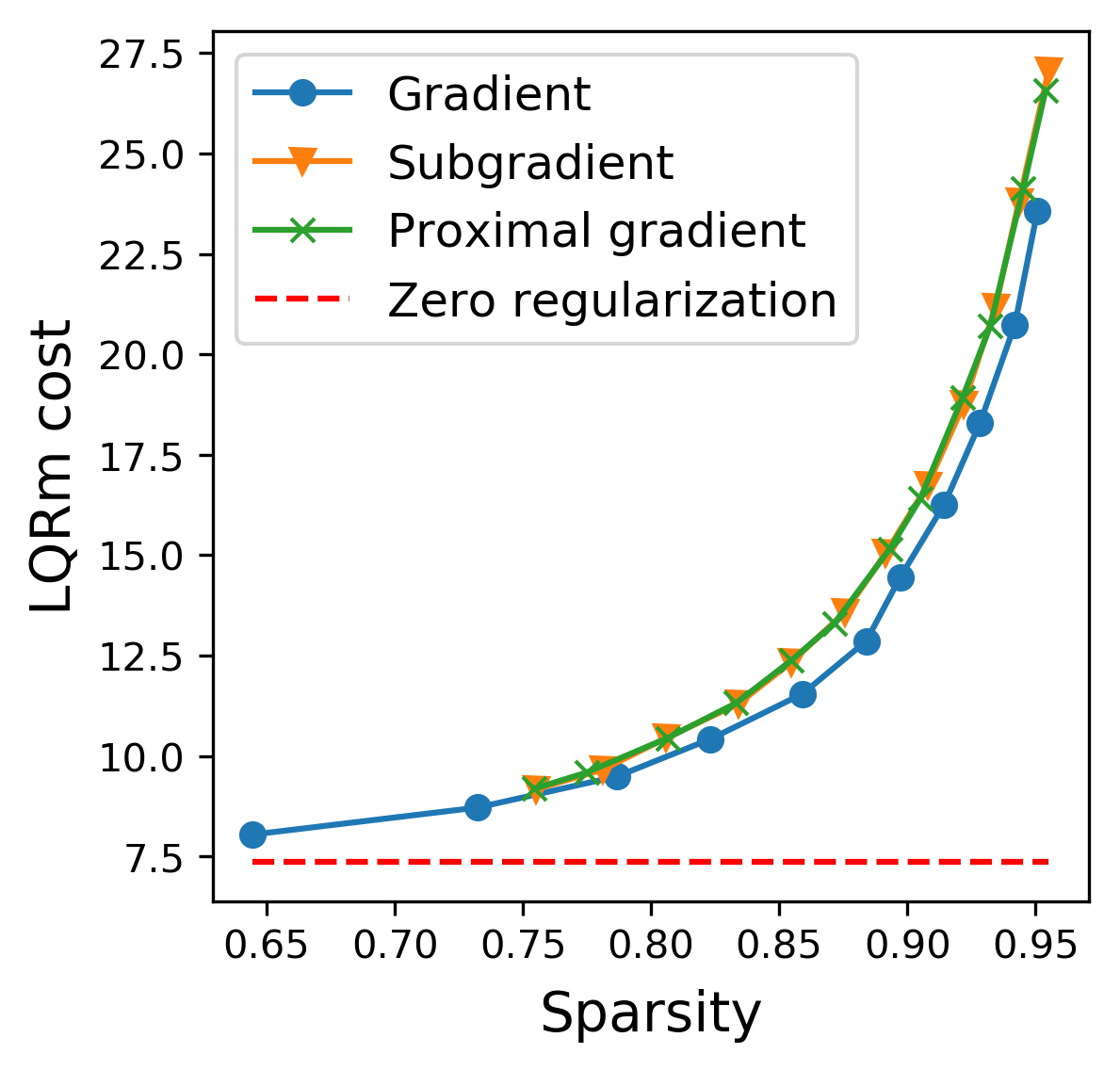}}
    \subfloat[High noise setting]{\includegraphics[width=4.3cm]{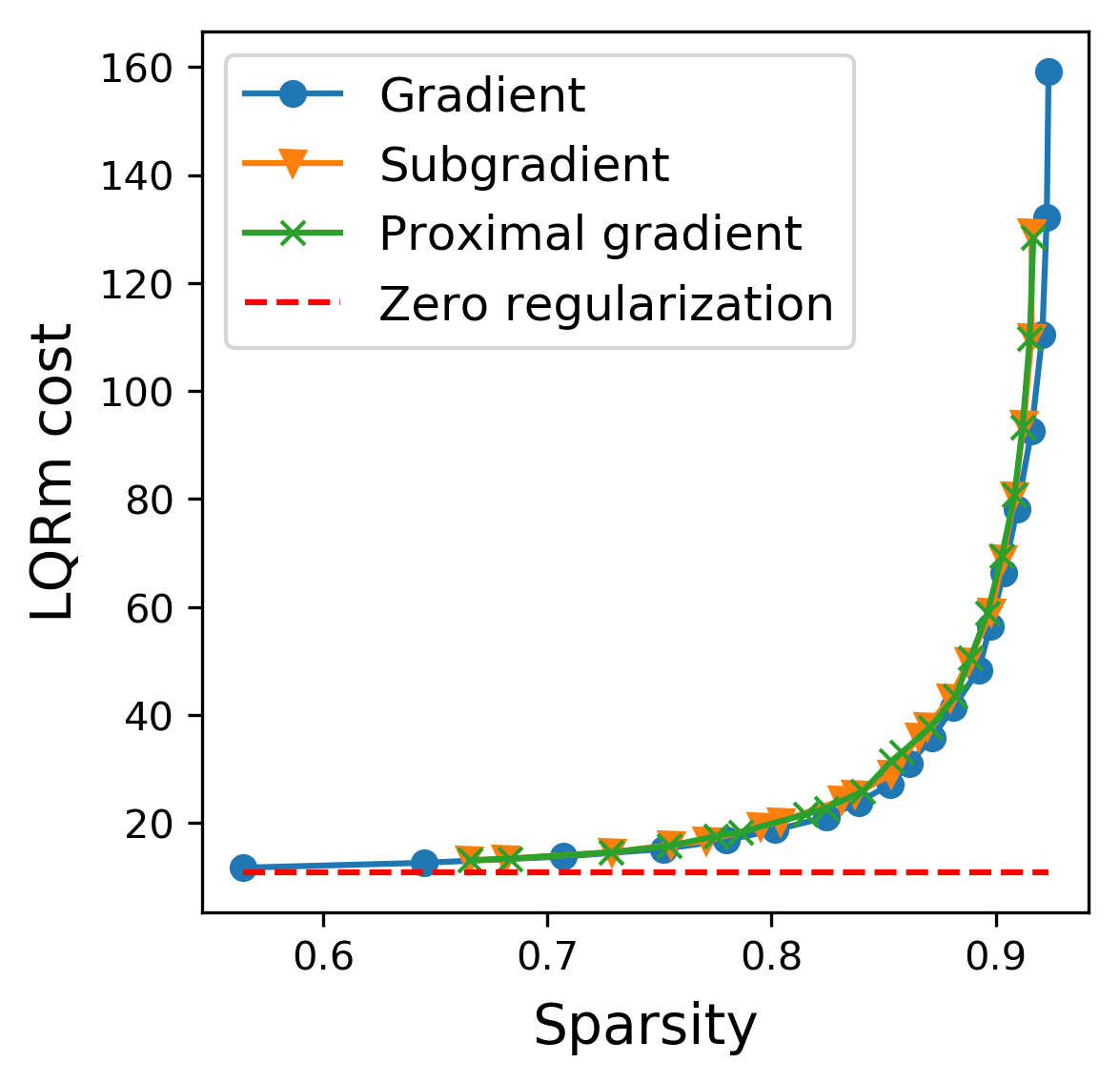}}
    \caption{LQRm cost vs. sparsity with $l_1$-norm regularization.}
    \label{fig:plot_comparison_vec1}
\end{figure}
\begin{figure}[ht!]
    \centering
    \subfloat[LQRm cost vs. sparsity]{\includegraphics[width=4.3cm]{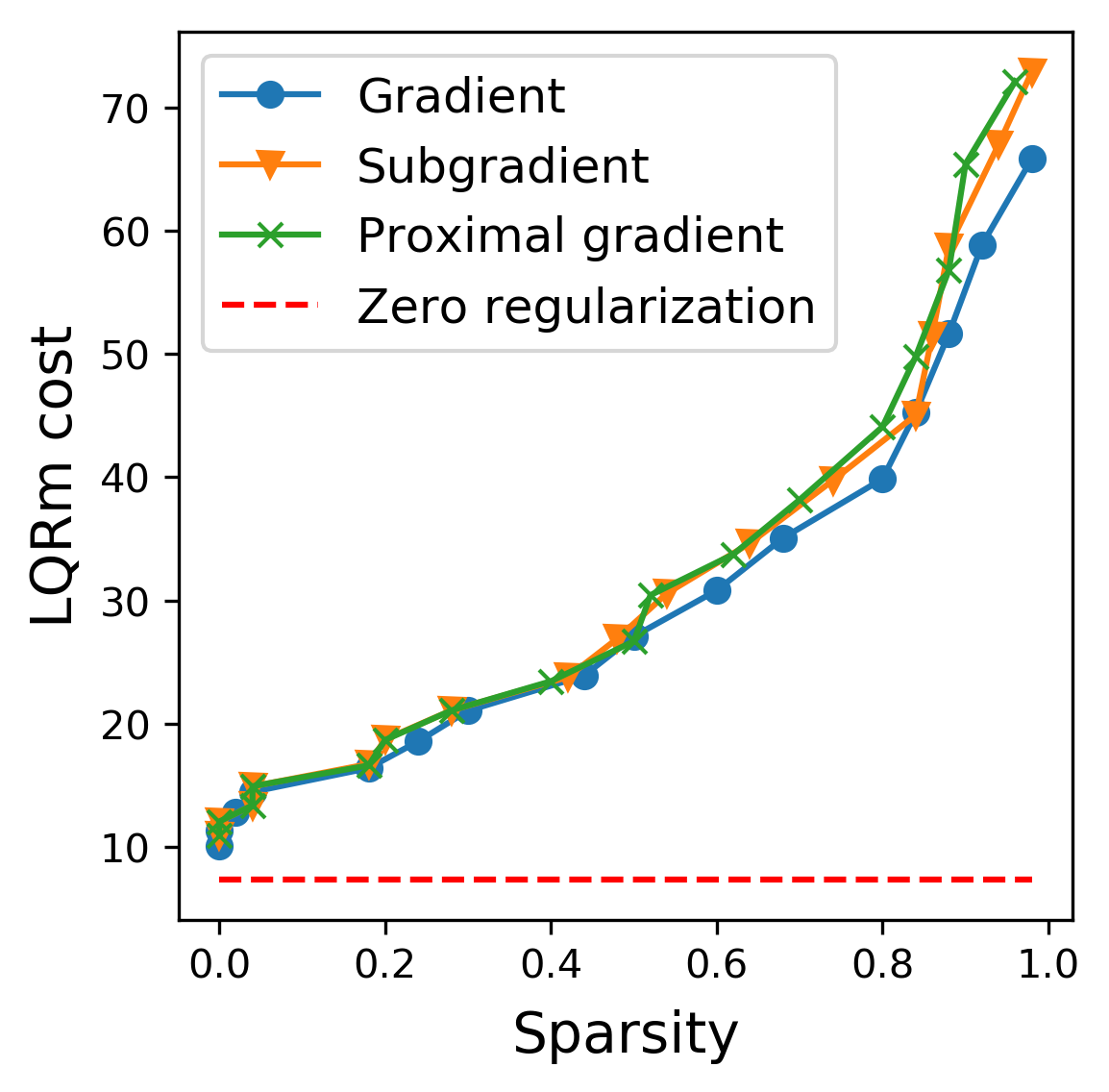}}
    \subfloat[Wall clock time vs. $\gamma$]{\includegraphics[width=4.3cm]{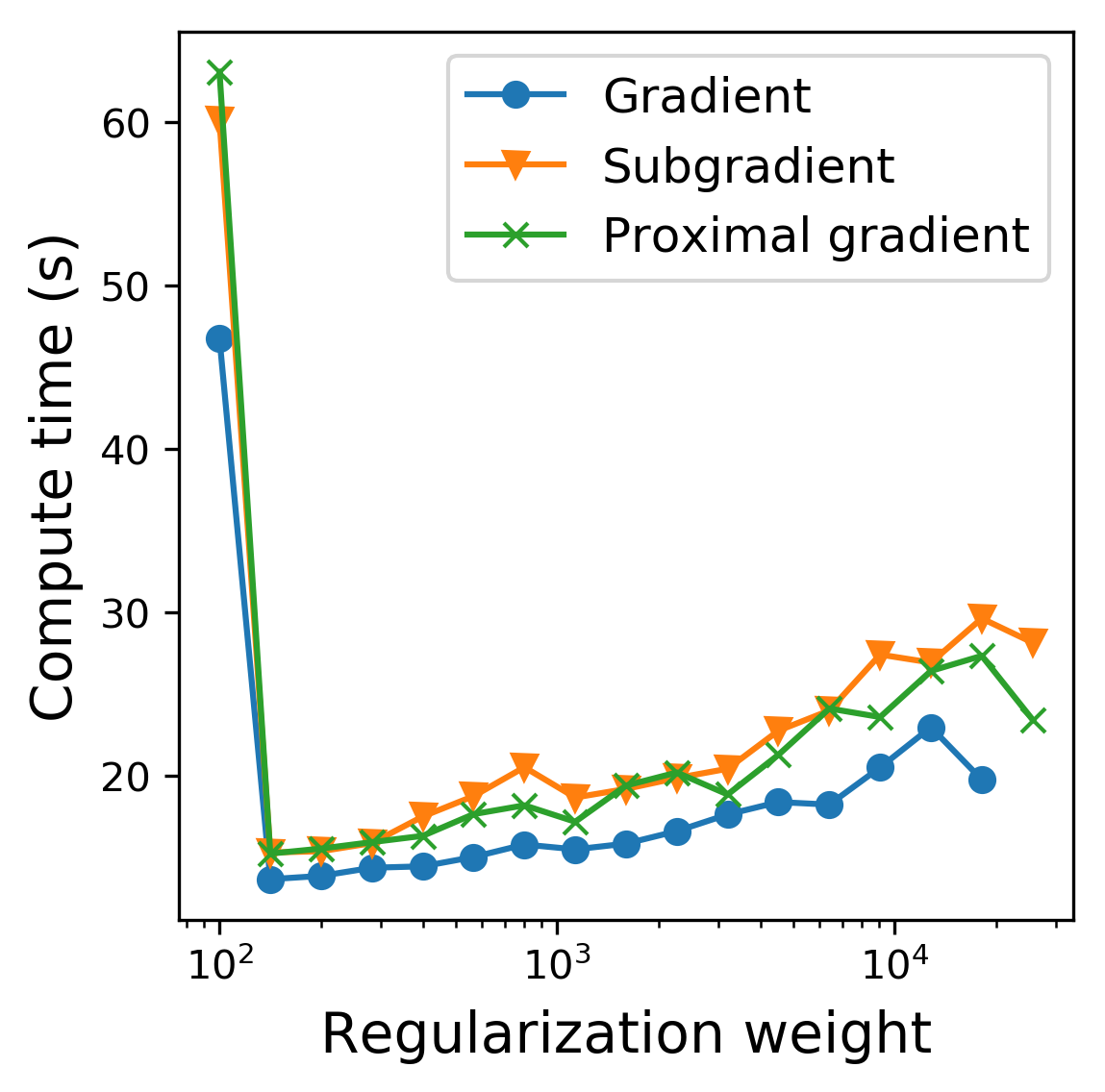}}
    \caption{Algorithm comparisons for the low noise setting with row group LASSO regularization.}
    \label{fig:plot_comparison_glr}
\end{figure}

For the subgradient and proximal gradient methods, we stopped iterating after the best iterate had been held for 100 iterations. For the gradient method, we stopped iterating when the Frobenius norm of the gradient of the cost function fell below a small threshold value, $0.1 \times \text{card}(K)$. We swept through a range of sparsity levels by solving a problem with low $\gamma$ then increasing $\gamma$ and resolving the problem using the previous solution as the initial guess. The step size $\eta$ was initialized at $10^{-5}$. For the $l_1$-norm and row group LASSO $\gamma$ was initialized at 10 and 100 respectively. For each successive problem, the regularization weight was multiplied by a ratio $r_\gamma=\sqrt{2}$ and the step size was multiplied by $r_\eta = r_\gamma^{-\sqrt[4]{2}}$.
To determine sparsity patterns we considered a value to be sparse if it was less than 5\% than the max value in $K$. For the $l_1$-norm the sparsity values were the absolute values of the entries. For the row group lasso norm the sparsity values were the the values are the $l_2$-norms of the rows and columns respectively.
Sparsity patterns are presented in Figs. \ref{fig:sparsity_vec1} and \ref{fig:sparsity_glr} with white cells representing near-zero entries.
The LQRm costs given in Figs. \ref{fig:plot_comparison_vec1} and \ref{fig:plot_comparison_glr} are for the sparse gains without any polishing step applied, which otherwise could significantly reduce the cost.
We give the total ``wall-clock'' computation time in Fig. \ref{fig:plot_comparison_glr} to capture the aggregate computational expense of each algorithm. 
The main computational expense came from evaluating the LQRm gradient at each iteration, which required solving a generalized discrete Lyapunov equation. 

As seen in Fig. \ref{fig:plot_comparison_glr}, the first iteration had the longest compute time since successive iterations benefited from favorable initial conditions from warm-starting. The compute time increased as the regularization weight was increased and a larger number of smaller steps were required to accommodate the increasing gradient magnitude.

From our empirical studies, the three methods presented all gave very similar results with similar efficacy; arbitrarily entrywise and row sparse mean-square stabilizing solutions were obtained for the low noise setting after a reasonable amount of computation time. Similarly, very sparse solutions for the high noise setting were obtained.

Python code which implements the algorithms and generates the figures reported in this work can be found in the GitHub repository at \url{https://github.com/TSummersLab/polgrad-multinoise/}.

The code was run on a desktop PC with a quad-core Intel i7 6700K 4.0GHz CPU, 16GB RAM.

\section{Concluding Remarks}
We developed three policy gradient algorithms for solving the sparse gain design problem for networked dynamical systems with multiplicative noise. 
We showed that the regularized LQR cost does not necessarily have a unique local minimum, hampering efforts to guarantee global convergence of the algorithms.
Nevertheless, efficacy of the algorithms is demonstrated empirically via computational simulations. Through various regularization functions we identified sparsity patterns for near-optimal actuator, sensor, and actuator-sensor link removal. This paves the way for data-driven control design in the model-free setting for such systems.

Future work will attempt to prove unique local minimization of the regularized LQR cost or provide a set of restrictions under which such a condition holds.
A salient issue with policy gradient methods relates to scalability; for large systems the gradient calculation is computationally expensive. Hence we will explore low-rank approximations of the gradient and consequent effects on convergence. We will also extend this work to the unknown-model setting and explore alternative model-based learning schemes.

\bibliography{ifacconf}             

\begin{thebibliography}{4}
\providecommand{\natexlab}[1]{#1}
\providecommand{\url}[1]{\texttt{#1}}
\providecommand{\urlprefix}{URL }
\expandafter\ifx\csname urlstyle\endcsname\relax
  \providecommand{\doi}[1]{doi:\discretionary{}{}{}#1}\else
  \providecommand{\doi}{doi:\discretionary{}{}{}\begingroup
  \urlstyle{rm}\Url}\fi

\bibitem[{Able(1956)}]{Abl:56}
Able, B. (1956).
\newblock Nucleic acid content of microscope.
\newblock \emph{Nature}, 135, 7--9.

\bibitem[{Able et~al.(1954)Able, Tagg, and Rush}]{AbTaRu:54}
Able, B., Tagg, R., and Rush, M. (1954).
\newblock Enzyme-catalyzed cellular transanimations.
\newblock In A.~Round (ed.), \emph{Advances in Enzymology}, volume~2, 125--247.
  Academic Press, New York, 3rd edition.

\bibitem[{Keohane(1958)}]{Keo:58}
Keohane, R. (1958).
\newblock \emph{Power and Interdependence: World Politics in Transitions}.
\newblock Little, Brown \& Co., Boston.

\bibitem[{Powers(1985)}]{Pow:85}
Powers, T. (1985).
\newblock Is there a way out?
\newblock \emph{Harpers}, 35--47.

\end{thebibliography}


\begin{thebibliography}{25}
\providecommand{\natexlab}[1]{#1}
\providecommand{\url}[1]{\texttt{#1}}
\providecommand{\urlprefix}{URL }
\expandafter\ifx\csname urlstyle\endcsname\relax
  \providecommand{\doi}[1]{doi:\discretionary{}{}{}#1}\else
  \providecommand{\doi}{doi:\discretionary{}{}{}\begingroup
  \urlstyle{rm}\Url}\fi

\bibitem[{Bollob{\'a}s and B{\'e}la(2001)}]{Bollobas2001}
Bollob{\'a}s, B. and B{\'e}la, B. (2001).
\newblock \emph{Random graphs}.
\newblock 73. Cambridge university press.

\bibitem[{Damm(2004)}]{Damm2004}
Damm, T. (2004).
\newblock \emph{Rational matrix equations in stochastic control}, volume 297.
\newblock Springer Science \& Business Media.

\bibitem[{Fazel et~al.(2018)Fazel, Ge, Kakade, and Mesbahi}]{Fazel2018}
Fazel, M., Ge, R., Kakade, S., and Mesbahi, M. (2018).
\newblock Global convergence of policy gradient methods for the linear
  quadratic regulator.
\newblock In \emph{Proceedings of the 35th International Conference on Machine
  Learning}, volume~80 of \emph{Proceedings of Machine Learning Research},
  1467--1476. PMLR.

\bibitem[{Gravell et~al.(2019)Gravell, Guo, and
  Summers}]{Gravell2019unpublished}
Gravell, B., Guo, Y., and Summers, T. (2019).
\newblock Policy gradient methods for networked dynamical systems with
  multiplicative noise.
\newblock
  \urlprefix\url{https://www.utdallas.edu/~tyler.summers/papers/Gravell2019TechReport.pdf}.
\newblock Unpublished.

\bibitem[{Hassan-Moghaddam and Jovanovi{\'c}(2018)}]{Hassan2018}
Hassan-Moghaddam, S. and Jovanovi{\'c}, M.R. (2018).
\newblock On the exponential convergence rate of proximal gradient flow
  algorithms.
\newblock In \emph{2018 IEEE Conference on Decision and Control (CDC)},
  4246--4251. IEEE.

\bibitem[{{Hassibi} et~al.(1998){Hassibi}, {How}, and {Boyd}}]{Hassibi1998}
{Hassibi}, A., {How}, J., and {Boyd}, S. (1998).
\newblock Low-authority controller design via convex optimization.
\newblock In \emph{Proceedings of the 37th IEEE Conference on Decision and
  Control (Cat. No.98CH36171)}, volume~1, 140--145 vol.1.

\bibitem[{Jadbabaie et~al.(2018)Jadbabaie, Olshevsky, and
  Siami}]{jadbabaie2018deterministic}
Jadbabaie, A., Olshevsky, A., and Siami, M. (2018).
\newblock Deterministic and randomized actuator scheduling with guaranteed
  performance bounds.
\newblock \emph{arXiv preprint arXiv:1805.00606}.

\bibitem[{Jovanovi{\'c} and Dhingra(2016)}]{jovanovic2016controller}
Jovanovi{\'c}, M.R. and Dhingra, N.K. (2016).
\newblock Controller architectures: Tradeoffs between performance and
  structure.
\newblock \emph{European Journal of Control}, 30, 76--91.

\bibitem[{Karimi et~al.(2016)Karimi, Nutini, and Schmidt}]{Karimi2016}
Karimi, H., Nutini, J., and Schmidt, M. (2016).
\newblock Linear convergence of gradient and proximal-gradient methods under
  the polyak-{\l}ojasiewicz condition.
\newblock In \emph{Machine Learning and Knowledge Discovery in Databases},
  795--811. Springer International Publishing, Cham.

\bibitem[{Kurdyka(1998)}]{Kurdyka1998}
Kurdyka, K. (1998).
\newblock On gradients of functions definable in o-minimal structures.
\newblock \emph{Annales de l'institut Fourier}, 48(3), 769--783.

\bibitem[{Lin et~al.(2013)Lin, Fardad, and Jovanovi{\'c}}]{Lin2013}
Lin, F., Fardad, M., and Jovanovi{\'c}, M.R. (2013).
\newblock Design of optimal sparse feedback gains via the alternating direction
  method of multipliers.
\newblock \emph{IEEE Transactions on Automatic Control}, 58(9), 2426--2431.

\bibitem[{Liu et~al.(2011)Liu, Slotine, and
  Barab{\'a}si}]{liu2011controllability}
Liu, Y.Y., Slotine, J.J., and Barab{\'a}si, A.L. (2011).
\newblock Controllability of complex networks.
\newblock \emph{Nature}, 473(7346), 167--173.

\bibitem[{Nesterov(2013)}]{Nesterov2013}
Nesterov, Y. (2013).
\newblock \emph{Introductory lectures on convex optimization: A basic course},
  volume~87.
\newblock Springer Science \& Business Media.

\bibitem[{Olshevsky(2014)}]{olshevsky2014minimal}
Olshevsky, A. (2014).
\newblock Minimal controllability problems.
\newblock \emph{IEEE Transactions on Control of Network Systems}, 1(3),
  249--258.

\bibitem[{Parikh et~al.(2014)Parikh, Boyd et~al.}]{Parikh2014}
Parikh, N., Boyd, S., et~al. (2014).
\newblock Proximal algorithms.
\newblock \emph{Foundations and Trends in Optimization}, 1(3), 127--239.

\bibitem[{Pasqualetti et~al.(2014)Pasqualetti, Zampieri, and
  Bullo}]{Pasqualetti2014c}
Pasqualetti, F., Zampieri, S., and Bullo, F. (2014).
\newblock Controllability metrics, limitations and algorithms for complex
  networks.
\newblock \emph{IEEE Transactions on Control of Network Systems}, 1(1), 40--52.

\bibitem[{Polyak et~al.(2013)Polyak, Khlebnikov, and
  Shcherbakov}]{Polyak-LMI_sparse_fb}
Polyak, B., Khlebnikov, M., and Shcherbakov, P. (2013).
\newblock An {LMI} approach to structured sparse feedback design in linear
  control systems.
\newblock In \emph{Proc. European Control Conference}, 833--838.

\bibitem[{Polyak(1963)}]{Polyak1963}
Polyak, B. (1963).
\newblock Gradient methods for the minimisation of functionals.
\newblock \emph{USSR Computational Mathematics and Mathematical Physics}, 3(4),
  864 -- 878.

\bibitem[{Ruths and Ruths(2014)}]{ruths2014control}
Ruths, J. and Ruths, D. (2014).
\newblock Control profiles of complex networks.
\newblock \emph{Science}, 343(6177), 1373--1376.

\bibitem[{Summers et~al.(2016)Summers, Cortesi, and
  Lygeros}]{summers2014submodularity}
Summers, T., Cortesi, F., and Lygeros, J. (2016).
\newblock On submodularity and controllability in complex dynamical networks.
\newblock \emph{IEEE Transactions on Control of Network Systems}, 3(1),
  91--101.

\bibitem[{Summers(2016)}]{summers2016actuator}
Summers, T. (2016).
\newblock Actuator placement in networks using optimal control performance
  metrics.
\newblock In \emph{IEEE Conference on Decision and Control}, 2703--2708.

\bibitem[{Taha et~al.(2019)Taha, Gatsis, Summers, and Nugroho}]{Taha2017d}
Taha, A.F., Gatsis, N., Summers, T., and Nugroho, S. (2019).
\newblock Time-varying sensor and actuator selection for uncertain
  cyber-physical systems.
\newblock \emph{IEEE Transactions on Control of Network Systems}.
\newblock \textit{to appear}.

\bibitem[{Tibshirani(1996)}]{Tibshirani1996}
Tibshirani, R. (1996).
\newblock Regression shrinkage and selection via the lasso.
\newblock \emph{Journal of the Royal Statistical Society. Series B
  (Methodological)}, 58(1), 267--288.

\bibitem[{Tzoumas et~al.(2016)Tzoumas, Rahimian, Pappas, and
  Jadbabaie}]{tzoumas2016}
Tzoumas, V., Rahimian, M.A., Pappas, G., and Jadbabaie, A. (2016).
\newblock Minimal actuator placement with bounds on control effort.
\newblock \emph{IEEE Transactions on Control of Network Systems}, 3(1), 67--78.

\bibitem[{Zare and Jovanovi{\'c}(2018)}]{Zare2018CDC}
Zare, A. and Jovanovi{\'c}, M.R. (2018).
\newblock Optimal sensor selection via proximal optimization algorithms.
\newblock In \emph{2018 IEEE Conference on Decision and Control (CDC)},
  6514--6518. IEEE.

\end{thebibliography}

\end{document}